\input louC.sty

\documentclass[a4paper,draft]{amsproc}
\usepackage{amsrefs}
\usepackage{amsmath}
\usepackage{mathrsfs}
\usepackage{amssymb}
\usepackage{amscd} 

\theoremstyle{plain}
 \newtheorem{thm}{\textbf{Theorem}}[section]

\theoremstyle{definition}

\theoremstyle{remark}
 
 \numberwithin{equation}{section}

\renewcommand{\leq}{\leqslant}
\renewcommand{\geq}{\geqslant}

\title[Asymptotic Bohr Radius for the Polynomials]{Asymptotic Bohr Radius for the Polynomials in One Complex Variable}

\subjclass[2010]{Primary 41}

\author[Chu]{\bfseries Cheng Chu}

\address{
Department of Mathematics \\ 
Washington University in Saint Louis  \\ 
Saint Louis, Missouri \\
USA}
\email{chengchu@math.wustl.edu}

\thanks{Partially supported by National Science Foundation Grant DMS 1300280}

\begin{document}

\vspace{18mm}
\setcounter{page}{1}
\thispagestyle{empty}

\begin{abstract}
We consider the Bohr radius $R_n$ for the class of complex polynomials in one variable of degree at most $n$. It was conjectured by R. Fournier in 2008 that $R_n={1\over 3}+{\pi^2\over {3n^2}}+o({1\over n^2})$. We shall prove this conjecture is true in this paper.
\end{abstract}

\maketitle

\section{Introduction}  
Let $\DD$ be the open unit disk in the complex plane $\CC$ and $H^\infty$ be the Banach space of bounded analytic functions on $\DD$ with the norm
$$
||f||_{\infty}=\sup_{z\in\DD}|f(z)|.
$$
Also let $\cP_n$ denote the subspace of $H^\infty$ consisting of  all the complex polynomials of degree at most $n$.
The Bohr radius $R$ for $H^\infty$ is defined as
$$
R=\sup\{r\in(0,1): \sum_{k=0}^\infty |a_k|r^k\leq ||f||_\infty,\, \m{for all}\, f(z)=\sum_{k=0}^\infty a_kz^k \in H^\infty\}.
$$
Bohr's famous power series theorem \cite{boh14} shows that $R={1\over 3}.$

In 2004, Guadarrama \cite{gua07} considered the Bohr type radius for the class $\cP_n$ defined by
\beq\label{def}
R_n=\sup\{r\in(0,1): \sum_{k=0}^n |a_k|r^k\leq ||p||_\infty,\, \m{for all}\, p(z)=\sum_{k=0}^n a_kz^k \in \cP_n\},
\eeq
and gave the estimate
$$
{{C_1}\over {3^{n/2}}}<R_n-{1\over 3}<C_2 {{\log n}\over n}, $$
for some positive constants $C_1$ and $C_2$. Later in 2008, Fournier obtained an explicit formula for $R_n$ by using the notion of bounded preserving operators. He proved the following theorem \cite{fou08}
\begin{thm}\label{fou}
For each $n \geq 1$, let $T_n(r)$ be the following $(n+1)\times(n+1)$ symmetric Toeplitz matrix
\beq\label{T}
\left(
  \begin{array}{cccccc}
    1 & r & -r^2 & r^3 & \cdots & (-1)^{n-1}r^n \\
    r & 1 & r & -r^2 & \cdots & (-1)^{n-2}r^{n-1} \\
    -r^2 & r & 1 & r &  &  \\
    r^3 & -r^2 & r & 1 & \ddots & \vdots \\
    \vdots &  & \ddots & \ddots & \ddots & \\
   (-1)^{n-1}r^n &  &\cdots  & & r & 1 \\
  \end{array}
\right).
\eeq
Then $R_n$ is equal to the smallest root in $(0,1)$ of the equation
$$
\det T_n(r)=0.
$$
\end{thm}
Based on the numerical evidence, he conjectured that
$$
R_n={1\over 3}+{\pi^2\over {3n^2}}+...
$$
The purpose of this note is to provide a positive answer. We shall prove
\begin{thm}\label{main}
Let $R_n$ be as in \eqref{def}, then
$$\lim_{n\to\infty} n^2\left(R_n-{1\over 3}\right)={\pi^2\over 3}.$$
\end{thm}
\section{Main Theorem}
In this section, we prove Theorem \ref{main}. The methods we use are similar to that in \cite{gre58}*{Chapter 5}.

{\bf Proof of Theorem \ref{main}.} Let $\Delta_n=\Delta_n(r)=\det T_n(r)$, where $T_n(r)$ is the symmetric Toeplitz matrix \eqref{T}. By Theorem \ref{fou}, $R_n$ is the smallest root in $(0,1)$ of the equation \beq\label{eq}\Delta_n(r)=0.\eeq

For $n\geq 2$, multiplying the second row of $\Delta_n$ by r, adding it to the first row and performing a similar operation with the columns, we have
\begin{align}
\nnb\Delta_n(r)=&\,\det
\left(
  \begin{array}{cccccc}
    1+3r^2 & 2r & 0 & & \cdots & 0 \\
    2r & 1 & r & -r^2 & \cdots & (-1)^{n-2}r^{n-1} \\
    0 & r & 1 & r &  &  \\
     & -r^2 & r & 1 & \ddots & \vdots \\
    \vdots &  & \ddots & \ddots & \ddots & \\
   0 &  &\cdots  & & r & 1 \\
  \end{array}
\right)\\\nnb\\
\label{re}=&\,(3r^2+1)\Delta_{n-1}(r)-4r^2\Delta_{n-2}(r).
\end{align}
If we set $\Delta_{-1}(r)=1$, then the recurrence relation \eqref{re} holds for all $n\geq 1$.

Consider the function associated with these Toeplitz matrices $\Delta_n$ $$f(r, \Gt)=1+\sum_{|n|>0}(-1)^{n-1}r^n e^{in\Gt}=\frac{3r^2+4r\cos\Gt+1}{r^2+2r\cos\Gt+1}.$$
In order to solve the equation \eqref{eq}, suppose
\beq\label{rep}
r=g(x)={1\over 3}(-2\cos x-\sqrt{4\cos^2 x-3}),
\eeq
for some $x\in [0, \pi]$. This substitution comes from $$3r^2+4r\cos x+1=0=f(r,x).$$
(In fact, for a fixed $r$, every eigenvalue $\Gl$ of $T_n$ can be written as $\Gl=f(r,x)$, for some $x\in [0,\pi]$. See \cite{gre58}*{Chapter 5}.)

Then \eqref{re} becomes
\beq\nnb
\Delta_n=(-4r\cos x)\Delta_{n-1}-4r^2\Delta_{n-2}.
\eeq
Its characteristic equation $$\Gl^2+4r\cos x\Gl+4r^2=0$$ has the roots $-2re^{\pm ix}$. Adding the initial conditions $\Delta_{-1}=\Delta_0=1$, we have
\beq\nnb
\Delta_n=\frac{(-2r)^{n+1}}{1-r^2}\left(\frac{\sin(n+2)x}{\sin x}+2r\frac{\sin(n+1)x}{\sin x}+r^2\frac{\sin nx}{\sin x}\right).
\eeq
Denote
\beq\nnb
p_n(\cos x)=\frac{\sin(n+2)x}{\sin x}+2r\frac{\sin(n+1)x}{\sin x}+r^2\frac{\sin nx}{\sin x}.
\eeq
Then $p_n(t)$ is a polynomial of degree $n+1$ in $t=\cos x$. Let $$x_\nu=\frac{\nu\pi}{n+2},\q\nu=1,2,\cdots,n+1.$$
Direct computation shows that $$p_n(\cos x_\nu)=(-1)^\nu 2r(1+r\cos\nu),$$ thus $$\m{sgn}\, p_n(\cos x_\nu)=(-1)^\nu.$$
Also $$\lim_{x\to 0^+}p_n(\cos x)>0.$$ So $p_n$ has $n+1$ distinct zeros $\{\cos t^{(n)}_\nu | \nu=1,2,\cdots, n+1\}$, such that
\beq\label{19}
0<t^{(n)}_1<x_1<t^{(n)}_2<x_2<\cdots<t^{(n)}_{n+1}<x_{n+1}<\pi.
\eeq
That means every root of the equation \eqref{eq} has the form \eqref{rep}. Notice that $g$ is positive only on $[{{5\pi}\over 6},\pi]$ and decreasing on $[{{5\pi}\over6},\pi]$, so the smallest root of \eqref{eq} in the interval $(0,1)$ is $g(t^{(n)}_{n+1})$, i.e. $R_n=g(t^{(n)}_{n+1})$.

Next, we will find an asymptotic expression for $t^{(n)}_{n+1}$. Notice that
\beq\label{20}
\lim_{n\to\infty}(-1)^{n+1}\frac{p_n(-\cos{z\over {n+2}})}{n+2}=(1-r)^2{{\sin z}\over z}.
\eeq
And \eqref{20} holds uniformly for $|z|<2\pi$. Let $$t^{(n)}_{n+1}=\frac{(n+1)\pi-\Ge_n}{n+2},$$ then $\Ge_n\in (0,\pi)$ by relation \eqref{19}.
Thus
\begin{align*}
0=&\lim_{n\to\infty}(-1)^{n+1}\frac{p_n(\cos{{(n+1)\pi-\Ge_n}\over {n+2}})}{n+2}\\
=&\lim_{n\to\infty}(-1)^{n+1}\frac{p_n(-\cos {\pi+\Ge_n\over n+2})}{n+2}\\
=&\lim_{n\to\infty}(1-r)^2{{\sin (\pi+\Ge_n)}\over \pi+\Ge_n}.
\end{align*}
Hence the accumulation point of $\{\Ge_n\}$ is either $0$ or $\pi$. Let$$y_n=\frac{(n+1)\pi-{\pi\over 2}}{n+2},\q y_n\in(x_n,x_{n+1}).$$
Using \eqref{20} again, we have
\begin{align*}
&\lim_{n\to\infty}(-1)^{n+1}\frac{p_n(\cos{{(n+1)\pi-{\pi\over 2}}\over {n+2}})}{n+2}\\
=&\lim_{n\to\infty}(-1)^{n+1}\frac{p_n(-\cos {{3\pi\over 2}\over n+2})}{n+2}\\
=&\lim_{n\to\infty}(1-r)^2{{\sin ({3\pi\over 2})}\over {3\pi\over 2}}<0.
\end{align*}
When $n$ is sufficiently large, $$\m{sgn}\,p_n(\cos y_n)=(-1)^n=\m{sgn}\,p_n(\cos x_n),$$
so $t^{(n)}_{n+1}\in (y_n, x_{n+1}).$
Consequently, $\Ge_n\to 0$ as $n\to\infty$, and then $$R_n=g(\pi-\Gt_n),$$ where
$$\Gt_n={{\pi+\Ge_n}\over n+2}={\pi\over n}+o({1\over n}),\q\m{as}\, n\to\infty.$$
By \eqref{rep},
\begin{align*}
R_n=&{1\over 3}(2\cos \Gt_n-\sqrt{4\cos^2 \Gt_n-3})\\
=&{1\over 3}+{\pi^2\over 3n^2}+o({1\over n^2}).
\end{align*}
\rightline{$\Box$}

\bibliography{references}

\begin{bibdiv}
\begin{biblist}

\bib{boh14}{article}{
      author={Bohr, H.},
       title={A theorem concerning power series},
        date={1914},
     journal={Proc. Lond. Math. Soc. (2)},
      volume={13},
       pages={1\ndash 5},
}

\bib{fou08}{article}{
      author={Fournier, R.},
       title={Asymptotics of the {Bohr} radius for polynomials of fixed
  degree},
        date={2008},
     journal={J. Math. Anal. Appl.},
      volume={338},
       pages={1100\ndash 1107},
}

\bib{gre58}{book}{
      author={Grenander, U.},
      author={Szeg{\"o}, G.},
       title={Toeplitz forms and their applications},
   publisher={University of California Press},
     address={Berkeley},
        date={1958},
}

\bib{gua07}{article}{
      author={Guadarrama, Z.},
       title={{Bohr}'s radius for polynomials in one complex variable},
        date={2007},
     journal={Comput. Methods Funct. Theory},
      volume={5},
       pages={143\ndash 151},
}

\end{biblist}
\end{bibdiv}

\end{document}